
\documentstyle{amsppt}
\baselineskip18pt
\magnification=\magstep1
\pagewidth{30pc}
\pageheight{45pc}
\hyphenation{co-deter-min-ant co-deter-min-ants pa-ra-met-rised
pre-print pro-pa-gat-ing pro-pa-gate
fel-low-ship Cox-et-er dis-trib-ut-ive}
\def\leaderfill{\leaders\hbox to 1em{\hss.\hss}\hfill}

\

\def\idest{i.e.,\ }

\def\a{{\alpha}}
\def\be{{\beta}}
\def\g{{\gamma}}

\def\bb{{\bold b}}

\def\bi{{\bold i}}
\def\bj{{\bold j}}

\def\bu{{\bold u}}

\def\brr{{\bar r}}
\def\b0{\text{\bf 0}}

\def\supp{\text{\rm \, supp}}

\def\boxit#1{\vbox{\hrule\hbox{\vrule \kern3pt
\vbox{\kern3pt\hbox{#1}\kern3pt}\kern3pt\vrule}\hrule}}
\def\rabbit{\vbox{\hbox{\kern0pt
\vbox{\kern0pt{\hbox{---}}\kern3.5pt}}}}

\def\tableau#1{
        \hbox {
                \hskip -10pt plus0pt minus0pt
                \raise\baselineskip\hbox{
                \offinterlineskip
                \hbox{#1}}
                \hskip0.25em
        }
}

\def\tabCol#1{
\hbox{\vtop{\hrule
\halign{\strut\vrule\hskip0.5em##\hskip0.5em\hfill\vrule\cr\lower0pt
\hbox\bgroup$#1$\egroup \cr}
\hrule
} } \hskip -10.5pt plus0pt minus0pt}

\def\CR{
        $\egroup\cr
        \noalign{\hrule}
        \lower0pt\hbox\bgroup$
}



\def\blank#1#2{
\hbox to #1{\hfill \vbox to #2{\vfill}}
}


\def\strut{\vrule height10pt depth5pt width0pt}

\topmatter
\title Schubert varieties and free braidedness
\endtitle

\author R.M. Green and J. Losonczy \endauthor
\affil Department of Mathematics \\ University of Colorado \\
Campus Box 395 \\ Boulder, CO  80309-0395 \\ USA \\ {\it  E-mail:}
rmg\@euclid.colorado.edu \\
\newline
Department of Mathematics\\ Long Island University\\
720 Northern Boulevard \\ Brookville,
NY  11548-1319\\ USA\\ {\it  E-mail:} losonczy\@liu.edu\\
\newline
\endaffil

\abstract We give a simple necessary and sufficient condition for
a Schubert variety $X_w$ to be smooth when $w$ is a freely braided
element of a simply laced Weyl group; such elements were
introduced by the authors in a previous work.  This generalizes in
one direction a result of Fan concerning varieties indexed by
short-braid avoiding elements. We also derive generating functions
for the freely braided elements that index smooth Schubert
varieties. All results are stated and proved only for the simply
laced case.
\endabstract


\endtopmatter

\centerline{\bf To appear in Transformation Groups}

\head Introduction \endhead

The ``freely braided elements'' of a simply laced Coxeter group
were defined in \cite{{\bf 7}}. The idea behind the definition is
that although it may be necessary to use long braid relations in
order to pass between two reduced expressions for a freely braided
element, the necessary long braid relations in a certain sense do
not interfere with one another.

In this paper, we study Schubert varieties indexed by freely
braided elements of a simply laced Weyl group; the
non-simply-laced case will not be considered.  Specifically, we
are interested in the problem of determining which of these
varieties are smooth. The general problem of characterising
smoothness for Schubert varieties has been investigated by
numerous authors in recent years. For example, Lakshmibai--Sandhya
\cite{{\bf 12}} have shown that in the type $A_n$ setting, a
Schubert variety $X_w$ is nonsingular if and only if the
permutation $w$ avoids the patterns $3412$ and $4231$. Billey
\cite{{\bf 1}} has obtained similar pattern-avoidance conditions
for rational smoothness in types $B_n$, $C_n$ and $D_n$.  General
necessary and sufficient conditions for smoothness have been given
by Kumar \cite{{\bf 11}}, and recently Billey--Postnikov
\cite{{\bf 3}} have found a criterion for smoothness in terms of
patterns in root systems. Important results have also been
obtained by Carrell, Peterson and others; for an overview of the
subject, the reader is referred to \cite{{\bf 2}}.

A main source of motivation for the present work is the paper
\cite{{\bf 6}}, in which Fan studies short-braid avoiding elements
of an arbitrary Weyl group and determines which of these elements
correspond to smooth Schubert varieties.  Fan found that for a
short-braid avoiding element $w$, the corresponding Schubert
variety $X_w$ is nonsingular if and only if $w$ is a product of
$\ell(w)$ distinct generators. Formulae for the enumeration of
such elements were also provided. Since in the simply laced
setting short-braid avoiding elements are a special case of freely
braided elements, it is natural to try to generalize Fan's results
to the freely braided case, and this is a principal goal of the
present article. Theorem 3.2 provides a particularly simple
characterisation of those freely braided $w$ for which $X_w$ is
smooth. A by-product of our argument is that we can describe
exactly when the deletion of a letter from a reduced expression
for a freely braided element results in another (not necessarily
freely braided) reduced expression (see Proposition 2.4 and Remark
2.5). The answer here is perhaps the simplest that can be
expected, and generalizes in one direction a deletion result of
Fan \cite{{\bf 6}, Theorem 1} concerning short-braid avoiding
elements of Weyl groups. Recently, Hagiwara--Ishikawa--Tagawa
\cite{{\bf 9}} have generalized the same deletion result in a
different way to a wider class of Coxeter groups.

In Theorem 4.2, we present generating functions for the number of
smooth varieties indexed by freely braided elements.  These may
easily be converted to explicit formulae using standard methods.

\head 1. Preliminaries \endhead

Let $W$ be a simply laced Coxeter group with set of distinguished
generators $S = \{s_i : i\in I\}$ and Coxeter matrix
$(m_{ij})_{i,j\in I}$.  Denote by $I^*$ the free monoid on $I$,
and let $\phi : I^* \longrightarrow W$ be the surjective morphism
of monoid structures satisfying $\phi(i) = s_i$ for all $i\in I$.
We say that a word $\bi \in I^*$ {\it represents} its image
$w=\phi(\bi)\in W$; further, if the length of $\bi$ (\idest the
number of factors used to express $\bi$ as a product of letters
from $I$) is minimal among the lengths of all the words that
represent $w$, then we call $\bi$ a {\it reduced expression} for
$w$. The length of $w$, denoted by $\ell(w)$, is then equal to the
length of $\bi$.

Let $V$ be a real vector space with basis $\{\a_i : i\in I\}$, and
denote by $B$ the {\it Coxeter form} on $V$ associated to $W$.
This is the symmetric bilinear form satisfying $B(\a_i,\a_j) = -
\cos{\frac{\pi}{m_{ij}}}$ for all $i,j\in I$. We view $V$ as the
underlying space of a reflection representation of $W$, determined
by the equalities $s_i \a_j = \a_j - 2B(\a_j,\a_i)\a_i$ for all
$i,j\in I$.

Denote by $\Phi$ the {\it root system} of $W$, \idest $\{w\a_i :
w\in W\text{ and }i\in I\}$.  Let $\Phi^+$ be the set of all $\be
\in \Phi$ such that $\be$ is expressible as a linear combination
of the $\a_i$ with nonnegative coefficients, and let $\Phi^- =
-\Phi^+$. We have $\Phi = \Phi^+ \cup \Phi^-$ (disjoint). The
elements of $\Phi^+$ (respectively, $\Phi^-$) are called {\it
positive} (respectively, {\it negative}) roots. The $\a_i$ are
also referred to as {\it simple} roots.

Associated to each $w\in W$ is the {\it inversion set} $\Phi(w) =
\Phi^+ \cap w^{-1}(\Phi^-)$, which has cardinality $\ell(w)$.
Given any reduced expression $i_1i_2\cdots i_n$ for $w$, we have
$\Phi(w) = \{r_1,r_2,\dots,r_n \}$, where $r_1 = \a_{i_n}$ and
$r_l = s_{i_n}\cdots s_{i_{n-l+2}}(\a_{i_{n-l+1}})$ for all $l\in
\{2,\dots,n\}$.  The sequence $\brr = (r_1,r_2,\dots,r_n)$ is
called the {\it root sequence} of $i_1i_2 \cdots i_n$, or a root
sequence {\it for} $w$. Note that any initial segment of a root
sequence is also a root sequence for some element of $W$.

Let $\bi, \bj \in I^*$ and let $i,j,k \in I$. We call the
substitution $\bi ij\bj \rightarrow \bi ji\bj$ a {\it commutation}
or {\it short braid move} if $m_{ij} = 2$, and if $m_{ij}=3$, we
call $\bi iji\bj \rightarrow \bi jij\bj$ a {\it long braid move}.
Applying a braid move to a reduced expression corresponds to
applying a permutation to the root sequence of that reduced
expression \cite{{\bf 7}, Proposition 3.1.1}. Specifically, let
$w\in W$ and suppose that $\bi i j \bj$ is a reduced expression
for $w$.  Let $\brr = (r_l)$ be the associated root sequence, and
let $n$ be the length of $\bj$. Then $m_{ij}=2$ if and only if
$r_{n+1}$ and $r_{n+2}$ are mutually orthogonal, in which case the
root sequence $\brr'$ of $\bi ji \bj$ can be obtained from $\brr$
by interchanging $r_{n+1}$ and $r_{n+2}$.  Employing again the
terminology used above for words, we say that the passage from
$\brr$ to $\brr'$ is obtained by a {\it commutation} or {\it short
braid move}. Suppose now that $\bi i j k \bj$ is a reduced
expression for $w$, and let $\brr = (r_l)$ be the associated root
sequence. Again, denote the length of $\bj$ by $n$.  Then $i = k$
if and only if $r_{n+1} + r_{n+3} = r_{n+2}$, in which case the
root sequence $\brr'$ of $\bi j i j \bj$ can be obtained from
$\brr$ by interchanging $r_{n+1}$ and $r_{n+3}$.  In this
instance, we say that the passage from $\brr$ to $\brr'$ is
obtained by a {\it long braid move}.

Let $w \in W$.  Any subset of $\Phi(w)$ of the form $\{ \a, \be,
\a + \be \}$ is called an {\it inversion triple} of $w$. We say
that an inversion triple $T$ of $w$  is {\it contractible} if
there is a root sequence for $w$ in which the elements of $T$
appear consecutively (in some order). The number of contractible
inversion triples of $w$ is denoted by $N(w)$.  If the
contractible inversion triples of $w$ are pairwise disjoint, then
$w$ is said to be {\it freely braided}.

\remark{Remark 1.1} Suppose that $W$ is of type $A_n$. Then the
freely braided elements $w \in W$ can be characterised by a
pattern avoidance condition involving four patterns (this was
pointed out by one of the referees, and is also discussed in
\cite{{\bf 7}, \S5.1}). More precisely, a permutation $w$ is
freely braided if and only if its 1-line notation avoids $3421$,
$4231$, $4312$ and $4321$. This fact will not be needed in the
sequel.
\endremark

\head 2. A deletion property of freely braided elements \endhead

We introduce a partial order $\leq$ on $\Phi$ by writing $\a \leq
\be$ if $\be - \a$ is a nonnegative linear combination of simple
roots. The following lemma will be used repeatedly in the proof of
Proposition 2.2.

\proclaim{Lemma 2.1} Let $w \in W$ and let $\brr$ be a root
sequence for $w$.  Suppose that $\brr = (\dots,\a,\be,\dots)$,
where $\a$ is not orthogonal to $\be$ relative to the Coxeter
form, and $\a \not\leq \be$. Then $\{\a - \be,\a,\be \}$ is a
contractible inversion triple of $w$.
\endproclaim

\demo{Proof} Let $\bi ji\bj$ be the reduced expression
corresponding to $\brr$, parsed in such a way that $\a =
\phi(\bj)^{-1}(\a_i)$ and $\be = \phi(\bj)^{-1}s_i(\a_j)$. We
first prove that $\{\a - \be,\a,\be \}$ is an inversion triple of
$w$ by showing that $\phi(\bj)$ has a reduced expression starting
with $j$. For this, it suffices by \cite{{\bf 10}, Proposition
5.7, Theorem 5.8} to verify that $\phi(\bj)^{-1}(\a_j)$ lies in
$\Phi^-$.

Since the root $\a$ is not orthogonal to $\be$, we have
$m_{ij}=3$, and hence $s_i(\a_j) = \a_j + \a_i$. We now compute
$\phi(\bj)^{-1}(\a_j) = \phi(\bj)^{-1}(\a_j + \a_i - \a_i) =
\phi(\bj)^{-1}s_i(\a_j) - \phi(\bj)^{-1}(\a_i) = \be-\a$, which
must lie in $\Phi^-$ since $\a \not\leq \be$.

Thus, $\{\a - \be,\a,\be\}$ is an inversion triple of $w$. Since
$\a$ and $\be$ are adjacent in $\brr$, \cite{{\bf 7}, Proposition
3.2.1} now implies that $\{\a - \be,\a,\be \}$ is contractible.
\qed\enddemo

We define the {\it height} of any root $\be$ to be the sum of the
coefficients used to express $\be$ as a linear combination of the
simple roots.

\proclaim{Proposition 2.2} Let $w \in W$ be freely braided and
suppose that $\a$ is the highest root of an inversion triple of
$w$.  Then $\a$ is the highest root of a contractible inversion
triple of $w$.
\endproclaim

\demo{Note} In the type $A_n$ setting, every inversion triple is
contractible, but this is not true in general (see \cite{{\bf 7},
Example 2.2.3, Proposition 5.1.1}).
\enddemo

\demo{Proof} Assume the contrary. Then, by Lemma 2.1, in every
root sequence for $w$, the root $\a$ cannot be directly to the
left of a root belonging to an inversion triple of $w$ in which
$\a$ is the highest root.  On the other hand, by hypothesis, $\a$
is the highest root of some inversion triple of $w$.  Choose a
root sequence $(\dots,\a,\be_1,\dots,\be_n,\g,\dots)$ for $w$ with
the property that $\{\a - \g,\a,\g\}$ is an inversion triple of
$w$, and $n \geq 1$ is as small as possible (over all inversion
triples of $w$ having $\a$ as highest root).

Now, $\be_n$ is not orthogonal to $\g$ by the minimality of $n$.
We claim that $\be_n \leq \g$.  Assume otherwise.  Then, according
to Lemma 2.1, $\{\be_n - \g,\be_n,\g\}$ is a contractible
inversion triple of $w$.  By \cite{{\bf 7}, Lemma 4.2.2}, we can
commute $\be_n - \g$ to the right, if necessary, so that it is
adjacent to $\be_n$. We cannot have $\be_n - \g = \a$, or else,
after the commutations, our sequence would take the form
$(\dots,\be_1,\dots,\be_{n-1},\a,\be_n,\g,\dots)$; a single long
braid move would then give us
$(\dots,\be_1,\dots,\be_{n-1},\g,\be_n,\a,\dots)$, with $\a$ lying
to the right of both $\g$ and $\a - \g$, an impossibility
\cite{{\bf 7}, Remark 2.2.2}.  So $\be_n - \g \neq \a$. But then
our sequence (after the commutations) looks like
$(\dots,\a,\be_1,\dots,\be_{n-1}, \be_n - \g,\be_n,\g,\dots)$,
where $\be_n - \g$ may or may not equal one of the $\be_l$. Either
way, performing a long braid move gives
$(\dots,\a,\be_1,\dots,\be_{n-1},\g,\be_n,\be_n - \g,\dots)$,
contradicting the minimality of $n$.  The claim is established,
\idest $\be_n \leq \g$.

Let $\g^*$ be a $\leq$-minimal element in $\{\be_l: \be_l \leq
\g\}$.  We claim that there is an index $m$ such that $\be_m$ is
to the left of $\g^*$ and $\be_m \not\perp \g^*$. If such an $m$
did not exist, then we could commute $\g^*$ to the left to obtain
the sequence $(\dots,\a,\g^*,\dots,\g,\dots)$. The root $\g^*$ is
not orthogonal to $\a$, because if it were, then we could commute
it past $\a$ and contradict the minimality of $n$. We also have
$\a \not\leq \g^*$ by the definition of $\g^*$. Hence, by Lemma
2.1, $\{\a - \g^*, \a, \g^*\}$ is a contractible triple,
contradicting the fact that $\a$ is not the highest root of any
contractible triple of $w$.  So there does indeed exist such an
$m$, and we may assume that $\be_m$ is directly to the left of
$\g^*$ in our chosen root sequence
$(\dots,\a,\be_1,\dots,\be_n,\g,\dots)$.

We have $\be_m \not\leq \g^*$ by our choice of $\g^*$. By Lemma
2.1 again, $\{\be_m - \g^*,\be_m,\g^*\}$ is a contractible triple
of $w$.  If necessary, we can commute $\be_m - \g^*$ to the right
in our chosen sequence so that it is adjacent to $\be_m$ (by
\cite{{\bf 7}, Lemma 4.2.2}). We can then apply a long braid move
so that the three roots of our triple appear in the order $\g^*,
\be_m, \be_m - \g^*$. If $\be_m - \g^*$ were to equal $\a$, then
this long braid move would place $\a$ closer to $\g$ than before,
a contradiction. So $\be_m - \g^* \neq \a$.

Suppose that $\be_m - \g^*$ was to the left of $\a$ before the
commutations and long braid move of the previous paragraph. Since
we commuted $\be_m - \g^*$ past $\a,\be_1, \dots,\be_{m-1}$ in
order to put it next to $\be_m$, it is orthogonal to all of these
vectors. We claim that the root $\g^*$ is also orthogonal to
$\a,\be_1, \dots,\be_{m-1}$. To see this, assume otherwise. Then,
after the long braid move of the previous paragraph, we can
commute $\g^*$ to the left until it is adjacent to a root $\delta$
among $\a,\be_1, \dots,\be_{m-1}$ with which it is not orthogonal.
Since $\g^*$ cannot be $\geq$ $\delta$, Lemma 2.1 implies that
$\g^*$ and $\delta$ belong to the same contractible triple of $w$.
This is a contradiction, since $\g^*$ cannot belong to any
contractible triple of $w$ other than $\{\be_m -
\g^*,\be_m,\g^*\}$, owing to the fact that $w$ is freely braided.
The claim is established. It follows that each root in $\{\be_m -
\g^*,\be_m,\g^*\}$ is orthogonal to $\a,\be_1, \dots,\be_{m-1}$,
and so can be commuted to the left past all of them, contradicting
the minimality of $n$.

Suppose finally that $\be_m - \g^*$ is one of the $\be_l$. Even
here, after applying the long braid move mentioned two paragraphs
above, we can commute $\g^*$ to the left past $\a$ (by the
argument of the previous paragraph), contradicting the minimality
of $n$. This last contradiction exhausts all possibilities for
$\be_m - \g^*$. \qed\enddemo

Let $\bi$ be a word in $I^*$ and suppose that $\bi$ can be written
as $\bu_0\bb_1 \bu_1\bb_2\bu_2 \cdots \bb_p\bu_p$, where each
$\bb_l$ is of the form $iji$ for some $i,j\in I$ with $m_{ij}=3$.
Then we call $\bb_1,\bb_2, \dots,\bb_p$ a {\it braid sequence} for
$\bi$. If $\bi \in I^*$ is reduced and $w=\phi(\bi)$ is freely
braided, then we say that $\bi$ is {\it contracted} provided there
exists a braid sequence for $\bi$ with $p=N(w)$ terms.

\proclaim{Proposition 2.3 \cite{{\bf 8}, Proposition 3.1.3}} Let
$w\in W$ be freely braided.
\item{\rm (i)}{There exists a contracted reduced expression
for $w$.}
\item{\rm (ii)}{Any contracted reduced expression for $w$ has
a unique braid sequence with $N(w)$ terms.} \qed\endproclaim

The proof of the next proposition adapts an argument of Fan
\cite{{\bf 6}, Theorem 1}.  We continue to assume that $W$ is
simply laced.

\proclaim{Proposition 2.4} Assume that $W$ is finite.  Let $w$ be
a freely braided element of $W$, and let $\bi= i_1\cdots i_n$ be a
contracted reduced expression for $w$ with braid sequence $\bb_1,
\dots,\bb_{N(w)}$. Then the expression obtained from $\bi$ by
deleting any letter, except a middle letter from one of the braids
$\bb_l$, is reduced.
\endproclaim

\demo{Proof} Suppose that the deletion of a letter $i_m$ from
$\bi$ results in a word $\bi'$ that is not reduced. We will show
that $i_m$ must be a middle letter of one of the braids $\bb_l$.
Denote by $y$ the element of $W$ represented by $i_{m+1} \cdots
i_n$, and note that $y^{-1}(\a_{i_m}) \in \Phi(w)$. Since $\bi$ is
contracted, we need only show that $y^{-1}(\a_{i_m})$ is the
highest root of a contractible triple of $w$.  By Proposition 2.2,
it suffices to show that $y^{-1}(\a_{i_m})$ is the highest root of
a (not necessarily contractible) inversion triple of $w$.

There is an index $p$ such that the expression $i_{p+1} \cdots
i_{m-1}i_{m+1} \cdots i_n$ is reduced and $i_p \cdots
i_{m-1}i_{m+1} \cdots i_n$ is not.  Let $x$ be the element of $W$
represented by the word $i_{p+1} \cdots i_{m-1}$, and let $\a$ be
the root $x^{-1}(\a_{i_p})$, which is positive since $i_p \cdots
i_{m-1}$ is reduced.

Let $\g_1 = y^{-1}(\a + \a_{i_m})$, $\g_2 = -y^{-1}(\a)$ and $\g_3
= y^{-1}(\a_{i_m})$.  Clearly, $\g_1 + \g_2 = \g_3$. We claim that
this is an inversion triple of the element $w'$ represented by
$i_p \cdots i_m \cdots i_n$.  Since $i_{p+1} \cdots i_{m-1}i_{m+1}
\cdots i_n$ is reduced and $i_p \cdots i_{m-1}i_{m+1} \cdots i_n$
is not, $y^{-1}x^{-1}(\a_{i_p})$ is negative.  This root equals
$y^{-1}(\a)$, hence $\g_2$ is positive.

Next, we show that $\g_1$ is positive.  Let $c = B(\a, \a_{i_m})$.
Note that $\a \neq \a_{i_m}$ (to see this, recall from above that
$y^{-1}(\a_{i_m})$ is positive whereas $y^{-1}(\a)$ is negative).
Therefore, since $W$ is finite, $c$ can only equal $-1/2$, $0$, or
$1/2$.  If $c$ were to equal $0$, then $y^{-1}s_{i_m}(\a)$ would
equal $y^{-1}(\a)$, a negative root; this would contradict the
fact that $y^{-1}s_{i_m}(\a)$ belongs to the inversion set of $w$.
If $c$ were to equal $1/2$, then the positive root
$y^{-1}s_{i_m}(\a)$ would equal $y^{-1}(\a - \a_{i_m}) =
y^{-1}(\a) - y^{-1}(\a_{i_m})$, which is negative, another
contradiction.  Thus, $c = -1/2$, and this implies that the
positive root $y^{-1}s_{i_m}(\a)$ equals $y^{-1}(\a + \a_{i_m}) =
\g_1$.

We verify that the roots $\g_1$ and $\g_2$ are sent negative by
$w'$. This follows from the calculations $w'(\g_1) =
s_{i_p}xs_{i_m}yy^{-1}(\a + \a_{i_m}) = s_{i_p}x(\a) = -\a_{i_p}$
and $ w'(\g_2) = -s_{i_p}xs_{i_m}yy^{-1}(\a) = -s_{i_p}x(\a +
\a_{i_m}) = -s_{i_p}(\a_{i_p} + x(\a_{i_m})) = \a_{i_p} -
s_{i_p}x(\a_{i_m})$. Regarding the latter, since $\a_{i_p}$ is
simple and $s_{i_p}x(\a_{i_m})$ is positive (because the word $i_p
\cdots i_m$ is reduced), we have that $w'(\g_2)$ is negative.

Thus, $\{\g_1,\g_2,\g_3\}$ is an inversion triple of $w'$. Since
every inversion triple of $w'$ is also an inversion triple of $w$,
this proves that $y^{-1}(\a_{i_m})$ is the highest root of some
inversion triple of $w$, as desired. \qed\enddemo

\remark{Remark 2.5} Let $w\in W$ be freely braided.  By \cite{{\bf
8}, Proposition 3.1.3(ii)}, any reduced expression for $w$ can be
transformed into a contracted expression by performing a sequence
of short braid moves. Thus, in the case where $W$ is finite,
Proposition 2.4 can actually be used to decide whether the
deletion of a letter from an arbitrary reduced expression for $w$
results in another reduced expression.
\endremark

\head 3. Application to Schubert varieties \endhead

Throughout this section, we assume that our simply laced group $W$
is the Weyl group of a semisimple simply-connected complex
algebraic group $G$ with fixed maximal torus $T$ and Borel
subgroup $B \supset T$. For each $w\in W$, let $X_w$ be the
associated Schubert variety, \idest the closure of the cell
$BwB/B$ in the generalized flag variety $G/B$. Define, for each
$w\in W$, the polynomial $P_w(t) = \sum_{v \preceq w}
t^{\ell(v)}$, where $\preceq$ denotes the Bruhat--Chevalley
partial order on $W$. Then $P_w(t^2)$ is the Poincar\'e polynomial
for the cohomology ring of $X_w$.

Carrell--Peterson have shown that a Schubert variety $X_w$ is
rationally smooth if and only if the polynomial $P_w(t)$ is
symmetric, meaning $P_w(t) = t^{\ell(w)} P_w(1/t)$ \cite{{\bf 4}}.
Note that in the simply laced setting, rational smoothness and
smoothness are equivalent for Schubert varieties (this is an
unpublished result of Peterson; see \cite{{\bf 5}}).

Let $w \in W$.  The set of letters from $I$ appearing in some
(any) reduced expression for $w$ will be called the {\it support}
of $w$, and will be denoted by $\supp(w)$. If $w$ is freely
braided, then it has a contracted reduced expression by
Proposition 2.3.  Consideration of this reduced expression reveals
that $\ell(w) - N(w) \geq \#\supp(w)$.

\definition{Definition 3.1}
If $w$ is freely braided and $\#\supp(w) = \ell(w) - N(w)$, then
$w$ is said to be {\it content maximal}.
\enddefinition

Let $v,w\in W$ and let $\bi = i_1 \cdots i_n$ be a reduced
expression for $w$. Recall that $v \preceq w$ if and only if $v$
is represented by a (possibly empty) word of the form
$i_{\mu_1}\cdots i_{\mu_m}$, where $1\leq \mu_1 < \cdots < \mu_m
\leq n$ \cite{{\bf 10}, Theorem 5.10}.

\proclaim{Theorem 3.2}  Let $w \in W$ be freely braided. Then the
Schubert variety $X_w$ is smooth if and only if $w$ is content
maximal.
\endproclaim

\demo{Proof}  By Proposition 2.4 and \cite{{\bf 6}, Lemma 2},
there are exactly $\ell(w) - N(w)$ elements of $W$ having length
$\ell(w) - 1$ that are less than $w$ in the Bruhat--Chevalley
order. If $X_w$ is (rationally) smooth, then $P_w(t)$ is
symmetric, and it follows that there are exactly $\ell(w) - N(w)$
elements of length $1$ that are less than $w$, as required.

Now assume that the support of $w$ has $\ell(w) - N(w)$ elements.

Let $\bi = \bu_0 \bb_1 \bu_1 \bb_2 \bu_2 \cdots \bb_{N(w)}
\bu_{N(w)}$ be a contracted reduced expression for $w$.  By
content maximality, the only letters appearing more than once in
$\bi$ are the outer factors of the braids $\bb_l$.  We will show
that the polynomial $P_w(t)$ is symmetric, and the theorem will
follow from the results of Carrell and Peterson cited above.

We argue by induction on $d = \ell(w)$.  For $d \leq 3$, the
polynomial $P_w(t)$ is symmetric by direct verification, so we
assume $d > 3$.  There are two cases to consider.

\medskip

\noindent {\bf Case 1:} $\bu_{N(w)}$ is a nonempty word, ending in
the letter $i$.

Let $\bi'$ be the reduced expression obtained from $\bi$ by
deleting $i$, and let $v$ be the group element represented by
$\bi'$. Now, $i$ appears exactly once in $\bi$, because it does
not appear in a braid. Therefore, by \cite{{\bf 8}, Lemma 3.2.2},
$v$ is freely braided and $N(v) = N(w)$.  We thus also have
$\#\supp(v) = \ell(v) - N(v)$, so that $v$ satisfies the inductive
hypothesis.

Since $i$ does not appear in $\bi'$, the set $\{x : x \preceq w\}$
equals the disjoint union of the sets $L_1 = \{x : x \preceq v\}$
and $L_2 = \{xs_i : x \preceq L_1\}$. Furthermore, for $x \in
L_1$, we have $\ell(xs_i) = \ell(x) + 1$. It follows that $P_w(t)
= (1+t) P_v(t)$. Since $\deg P_v(t) = \deg P_w(t) - 1$ and
$P_v(t)$ is symmetric by induction, $P_w(t)$ is symmetric.

\medskip

\noindent {\bf Case 2:} $\bu_{N(w)}$ is empty.

Since $d>3$, we have $N(w)\geq 1$. Write $\bb_{N(w)} = iji$ and
let $v = ws_is_js_i$. Since $i$ and $j$ do not lie in the support
of $v$, any element $x \preceq w$ is uniquely of the form $x = ab$
where $a \preceq v$, $b \preceq s_is_js_i$ and $\ell(ab) = \ell(a)
+ \ell(b)$.

It follows that $P_w(t) = (1 + t)(1 + t + t^2) P_v(t)$.  By
\cite{{\bf 8}, Lemma 3.2.2} again (applied three times
successively to $w$), $v$ is freely braided and $N(v) = N(w) - 1$.
Since $\ell(v) = \ell(w) - 3$ and $\#\supp(v) = \#\supp(w) - 2$,
the inductive hypothesis is satisfied by $v$, so that $P_v(t)$ is
symmetric. Since $\deg P_v(t) = \deg P_w(t) - 3$, the induction is
complete. \qed\enddemo

\remark{Remark 3.3} Another possible approach to proving Theorem
3.2 would be to use the criterion for smoothness recently found by
Billey--Postnikov \cite{{\bf 3}}.
\endremark

\medskip

An element $w\in W$ is said to be {\it fully commutative} if every
reduced expression for $w$ can be transformed into any other by
performing a sequence of short braid moves \cite{{\bf 14}}.  This
is equivalent to the requirement that $N(w) = 0$ \cite{{\bf 8},
Proposition 1.2.2}. Thus, every fully commutative element of $W$
is freely braided.  Note that the fully commutative elements
coincide with the ``short-braid avoiding" elements of \cite{{\bf
6}} in the case where $W$ is simply laced (the only case
considered in this paper).

We have the following immediate corollary of Theorem 3.2, which
was originally proved by Fan for short-braid avoiding elements of
an arbitrary Weyl group \cite{{\bf 6}, Proposition 3}.

\proclaim{Corollary 3.4}  Let $w \in W$ be fully commutative. Then
the Schubert variety $X_w$ is smooth if and only if $w$ equals a
product of $\ell(w)$ distinct generators. \qed\endproclaim

\head 4. Enumeration of content maximal elements \endhead

Motivated by Theorem 3.2, we derive generating functions for the
number of content maximal elements in types $A_n$, $D_n$ and
$E_n$. In the $E$-series, we allow $n$ to be any positive integer
$\geq 6$, since all groups in this series have finitely many
freely braided elements \cite{{\bf 8}, Theorem 3.3.3} and our
methods apply equally well to them.

It turns out that all three generating functions are based on a
single recurrence relation; we explain by using the following
set-up.  Consider a nested sequence of simply laced Coxeter groups
$ W_1 \subset W_2 \subset \cdots $ with generating sets $S_1
\subset S_2 \subset \cdots $ and alphabets $I_1 \subset I_2
\subset \cdots$\,. Thus, for each $n\geq 1$, we have $S_n = \{s_i
: i\in I_n \}$. Assume that for any positive integers $m < n$ and
any $i,j \in I_m$, the order of $s_is_j$ is the same in $W_m$ as
it is in $W_n$.  Assume further the existence of a sequence
$(i_n)_{n\geq 1}$ such that each $i_n$ lies in $I_n$, and whenever
$n > 1$, we have $m_{i_nj}=2$ for all $j\in I_{n-1}$ except for $j
= i_{n-1}$.

Some ad hoc terminology will be useful.  Given a content maximal
element $w$ and a letter $i$, we say that $i$ is {\it not braided}
in $w$ if $i$ appears at most once in any reduced expression for
$w$. If $i$ appears more than once in some reduced expression for
$w$, then by content maximality it appears exactly twice and there
is precisely one letter $j$ between the two occurrences of $i$
such that $m_{ij}=3$.  In this case, we say that $i$ is {\it
braided} in $w$ {\it with} $j$.  It is clear that if $i$ is
braided (in $w$) with $j$, then $j$ is braided with $i$. Note that
if $i$ is braided in $w$, then it is braided with precisely one
other letter.

Let $w\in W$.  In the proof of the following lemma, we will use a
result of Matsumoto \cite{{\bf 13}} and Tits \cite{{\bf 15}},
which states that any reduced expression for $w$ can be
transformed into any other by applying a sequence of long and
short braid moves.

Denote by $f(n)$ the number of content maximal elements in $W_n$.

\proclaim{Lemma 4.1} Let $w \in W_n$ be content maximal.  Then $w$
satisfies  one of the following seven mutually exclusive
conditions:
\item{\rm (i)}{$i_n \not\in \supp(w)${\rm ;}}
\item{\rm (ii)}{$i_n \in \supp(w)$ and $i_n$ is not braided in $w$ and
if $n>1$, then $i_{n-1} \not\in \supp(w)${\rm ;}}
\item{\rm (iii)}{$i_n \in \supp(w)$ and $i_n$ is not braided in $w$ and
$i_{n-1} \in \supp(w)$ and $\ell(ws_{i_n}) < \ell(w)${\rm ;}}
\item{\rm (iv)}{$i_n \in \supp(w)$ and $i_n$ is not braided in $w$ and
$i_{n-1} \in \supp(w)$ and $\ell(s_{i_n}w) < \ell(w)${\rm ;}}
\item{\rm (v)}{$i_n \in \supp(w)$ and $i_n$ is braided in $w$ and
if $n>2$, then $i_{n-2} \not\in \supp(w)${\rm ;}}
\item{\rm (vi)}{$i_n \in \supp(w)$ and $i_n$ is braided in $w$ and
$i_{n-2} \in \supp(w)$ and $\ell(ws_{i_{n-1}}) < \ell(w)${\rm ;}}
\item{\rm (vii)}{$i_n \in \supp(w)$ and $i_n$ is braided in $w$ and
$i_{n-2} \in \supp(w)$ and $\ell(s_{i_{n-1}}w) < \ell(w)$.}

For $n>3$, the number of content maximal elements of $W_n$
satisfying the respective conditions {\rm (i),} {\rm (ii), }
\dots, {\rm (vii)} is $f(n-1)$, $f(n-2)$, $f(n-1) - f(n-2)$,
$f(n-1) - f(n-2)$, $f(n-3)$, $f(n-2) - f(n-3)$, and $f(n-2) -
f(n-3)$. We therefore have the recurrence $f(n) = 3f(n-1) + f(n-2)
- f(n-3)$.
\endproclaim

\demo{Note} In conditions (iii)--(v), we are assuming $n
> 1$.  In conditions (vi) and (vii), we are assuming $n > 2$.
\enddemo

\demo{Proof} Suppose that $i_n \in \supp(w)$ and $i_n$ is not
braided in $w$ and $i_{n-1} \in \supp(w)$.  We claim that either
$\ell(ws_{i_n}) < \ell(w)$ or $\ell(s_{i_n}w) < \ell(w)$. Fix a
contracted reduced expression for $w$.  By our assumptions, $i_n$
appears exactly once and $i_{n-1}$ appears at least once and at
most twice in this expression. If $i_{n-1}$ appears only once,
then, since $m_{i_nj} = 2$ for all $j\in I_{n-1}$ except for
$j=i_{n-1}$, we can commute $i_n$ to one of the ends of our
contracted expression. If instead $i_{n-1}$ appears twice, then it
is braided, necessarily with $i_{n-2}$. Once again, we can commute
$i_n$ to an end of our contracted expression. Thus, either (iii)
or (iv) holds.  Note that these conditions cannot hold
simultaneously, by our assumption that $i_n$ is not braided in $w$
together with the result of Matsumoto and Tits cited above.

Similar reasoning shows that if $i_n \in \supp(w)$ and $i_n$ is
braided in $w$ and $i_{n-2} \in \supp(w)$, then exactly one of the
conditions (vi) and (vii) holds.  Thus, our seven conditions are
mutually exclusive and collectively exhaustive.

Regarding the cardinality assertions, let $C_n$ denote the set of
content maximal elements in $W_n$ and let $C_n$(i), $C_n$(ii),
\dots , $C_n$(vii) denote the respective subsets of elements
satisfying conditions (i), (ii), \dots , (vii).  We have $C_n$(i)
$\, = C_{n-1}$, so that $\# C_n$(i)$\, = f(n-1)$. The mapping
$C_n$(ii)$\, \longrightarrow C_{n-2}$ given by $w \mapsto
ws_{i_n}$ is a well-defined bijection, hence $\# C_n$(ii)
$=f(n-2)$. Similarly, there are bijections $C_n$(iii)$ \,
\longrightarrow C_{n-1}\setminus C_{n-2}$ and $C_n$(iv)$ \,
\longrightarrow C_{n-1}\setminus C_{n-2}$, given by $w \mapsto
ws_{i_n}$ and $w \mapsto s_{i_n}w$, respectively. It follows that
$\# C_n$(iii)$\, = \# C_n$(iv)$\, = f(n-1) - f(n-2)$. The
remaining bijections are induced by left or right multiplication
by $s_{i_n}s_{i_{n-1}}s_{i_n}$. \qed\enddemo

Through direct computation, using either a variation of the proof
of the above lemma or distinguished coset representatives for the
subgroup $W_{n-1}$ of $W_n$, one finds that the respective numbers
of content maximal elements in Coxeter groups of type $A_1$,
$A_2$, $A_3$, $D_4$, $D_5$ and $E_6$ are $2$, $6$, $19$, $62$,
$201$ and $652$. This information, together with Lemma 4.1 itself,
enables one to derive the ordinary generating functions for the
numbers $f(n)$ when $W_n$ is of type $A_n$ or $D_n$ or $E_n$, as
follows.

\proclaim{Theorem 4.2} The ordinary generating function for the
number of content maximal elements is
$$\cases \displaystyle{\frac{2x-x^3}{1-3x-x^2+x^3} = 2x + 6x^2
+ 19x^3 + 61x^4 + 196x^5 + \cdots } & \text{ in type }A_n, \cr
\displaystyle{\frac{62x^4 + 15x^5 - 19x^6}{1-3x-x^2+x^3} =
62x^4+201x^5+ 646x^6 + 2077x^7 + \cdots } & \text{ in type }D_n,
\cr \displaystyle{\frac{652x^6 + 140x^7 - 201x^8}{1-3x-x^2+x^3} =
652x^6 + 2096x^7 + 6739x^8 + \cdots } & \text{ in type }E_n.
\endcases$$ \qed\endproclaim

\head Acknowledgment \endhead

The authors thank the referees for their helpful and stimulating
comments.

\leftheadtext{} \rightheadtext{}

\end